\documentclass[11pt]{amsart}
\pdfoutput=1
\newtheorem{thm}{Theorem}[section]

\newtheorem{lem}[thm]{Lemma}
\newtheorem{cor}[thm]{Corollary}

\theoremstyle{definition}
\newtheorem{defn}[thm]{Definition}

\newtheorem{ques}[thm]{Question}

\usepackage[cp1250]{inputenc}
\usepackage{enumitem, graphicx, lpic}
\usepackage[sc]{mathpazo}
\usepackage{amsmath, amssymb, amsthm}

\begin{document}


\title{Knotted surfaces and equivalencies of their diagrams without triple points}

\author{Micha\L\ Jab\L onowski}
\date{\today}
\address{Institute of Mathematics, University of Gda\'nsk\\
Wita Stwosza 57, 80-952 Gda\'nsk, Poland\\
E-mail: michal.jablonowski@gmail.com}

\begin{abstract}
The singularity set of a generic standard projection to the three space of a closed surface linked in four space, consists of at most three types: double points, triple points or branch points. We say that this generic projection image is p-diagram if it does not contain any triple point. Two p-diagrams of equivalent surface links are called p-equivalent if there exists a finite sequence of local moves, such that each of them is one of the four moves taken from the seven on the well known Roseman list, that connects only p-diagrams. It is natural to ask: are any of two p-diagrams of equivalent surface links always p-equivalent? We introduce an invariant of p-equivalent diagrams and an example of linked surfaces that answers our question negatively.
\end{abstract}

\maketitle

\section{Introduction}

Let $F$ denotes a closed surface, the image of an embedding $f:F\to \mathbb{R}^4$ is called the \emph{knotted surface}, with $F$ being the \emph{based surface} for $f$.  Sometimes to emphesize that the based surface may not be connected, instead of a knotted surface we say a \emph{linked surface}.  For simplicity it is often used $F$ to denote $f(F)$. 

Two knotted surfaces are \emph{equivalent} if there exists an orientation preserving autohomeomorphism of $\mathbb{R}^4$, taking one surface to the other.

\begin{defn}
We say that a diagram $D_F$ is a \emph{$p$-diagram} if it does not contain a triple point, that is for all $x\in D_F$ we have $\#\left(\pi^{-1}(x)\cap F\right)<3$.  Two $p$-diagrams $D_F$ and $D_{F'}$ are \emph{$p$-equivalent} if there exists a finite sequence of $p$-diagrams 
$D_F=D_1\to D_2\to\cdots\to D_n=D_{F'}$ such that for all $i=1, 2, \ldots, n-1$ a transition from $D_i$ to $D_{i+1}$ is one of four Roseman moves that involve only $p$-diagrams (or some isotopy of that diagram in $\mathbb{R}^3$).
\end{defn}

Above mention four moves can be (of course) realized by an isotopy of a surface in $\mathbb{R}^4$ without introducing a triple point in the projection. It is natural then to consider the following problem.

\begin{ques}
Are any of two $p$-diagrams of equivalent linked surfaces always $p$-equivalent?
\end{ques}

Similar question but for branch points instead of triple points in our definition of $p$-diagram, was answered by S. Satoh in \cite{Sat01b}. In this paper we introduce an invariant of $p$-equivalent diagrams, that helps to answer our question as follows.

\begin{thm}\label{t2}
There exist two $p$-diagrams of equivalent surfaces that are not $p$-equivalent.
\end{thm}

\section{Preliminaries}

Let us fix a projection $\pi:\mathbb{R}^4\to\mathbb{R}^3$ by \emph{$\pi((x, y, z, t))=(x, y, z)$}. A \emph{double decker set} $\Gamma$ is the closure in $F$ of a set $\{p\in F:\#(\pi^{-1}(\pi(p))\cap F)>1\}$, and the \emph{double point set} is the image $\pi(\Gamma)$ which we will denote by $\Gamma^*$. We say that the $\pi(F)$ also denote by $F^*$ is in a \emph{general position} if the double point set consists of points whose neighborhood is locally homeomorphic to:

\begin{enumerate}[label=(\roman{*})]
\item two transversely intersecting sheets,
\item three transversely intersecting sheets,
\item the Whitney'a umbrella.
\end{enumerate}

\begin{figure}[ht]
\begin{center}
\begin{lpic}[]{r31(11cm)}
   \lbl[t]{11,0;(i)}
   \lbl[b]{43,20;(ii)}
   \lbl[t]{78,0;(iii)}
\end{lpic}
\caption{\label{r31}}
\end{center}
\end{figure}
Points corresponding to these cases (i), (ii), (iii) in the Figure \ref{r31} are called a \emph{double point}, a \emph{triple point} and a \emph{branch point} respectively.

\begin{thm}[\cite{Ros00}]
Every knotted surface is equivalent (and is placed arbitrarily near) to a knotted surface whose projection is in a general position.
\end{thm}

We can therefore from now on assume a general position of a projection image.

A \emph{diagram} $D_K$ of a knotted surface $K$ is the image $\pi(K)$ with additional information at self crossing, about which sheet was higher before being projected.

We have the following well known characterization of the equivalency between knotted surfaces. 

\begin{thm}[\cite{Ros98}]
Two diagrams represent equivalent linked surfaces if and only if one of them may be achieved from the other  by a finite sequence of elementary local moves taken from the list in the Figure \ref{r43} (and some isotopy of the diagram in $\mathbb{R}^3$).
\end{thm}

\begin{figure}[ht]
\begin{center}
\begin{lpic}[]{r43(11cm)}
   \lbl[l]{25,100;\emph{I}}
   \lbl[l]{55,101;\emph{I}\emph{I}}
   \lbl[l]{20,66;\emph{I}\emph{I}\emph{I}}
   \lbl[l]{62,65;\emph{I}\emph{V}}
   \lbl[l]{94,80;V}
   \lbl[l]{25,24;VI}
   \lbl[b]{70,30;VII}
\end{lpic}
\caption{\label{r43} (see \cite{CarSai98})}
\end{center}
\end{figure}

For the simplicity, in pictures of those \emph{Rosemann moves} we do not indicate which sheet is an upper sheet.

\begin{defn}
Let $F\subset\mathbb{R}^4$ be an oriented linked surface such that the diagram $D_F$ is a $p$-diagram.  Define \emph{$A_f$} as a set of those circles $c\subset\Gamma$, that the set $\left(\pi^{-1}\left(\pi(c)\right)\right)\cap F$  has nonempty intersection with two distinct components of the surface $F$.
\end{defn}

On each circle in the set $A_f$ let us give an orientation defined as follows and depicted in Figure \ref{moje3}.

\begin{defn}[\cite{Sat01b}]
Let $G_1^*$ and $G_2^*$ be sheets in $F^*$, intersecting each other in a double point curve $c^*$. Let $c_1\subset G_1$ and $c_2\subset G_2$ be curves on the double decker set, that are projection preimages of the curve $c^*\subset\Gamma^* $. Let $\vec{n_1}$ and $\vec{n_2}$ be normal vectors to $G_1^*$ and $G_2^*$ respectively. Then the \emph{orientation of the vector $\vec{v_1}$} on $c_1$ is defined through the condition, that a triple $(\vec{n_1},\vec{n_2},\pi(\vec{v_1}))$ matches the orientation given to $\mathbb{R}^3$. By analogy we can define the orientation of $\vec{v_2}$ on $c_2$.
\end{defn}
\begin{figure}[ht]
\begin{center}
\begin{lpic}[]{moje3(11cm)}
   \lbl[l]{3,20;$F^*$}
   \lbl[l]{73,20;$F$}
   \lbl[r]{5,40;$G_1^*$}
   \lbl[l]{-5,30;$G_2^*$}
   \lbl[t]{44,27;$G_1$}
   \lbl[t]{68,29;$G_2$}
   \lbl[l]{52,34;$c_1$}
   \lbl[l]{76,33;$c_2$}
   \lbl[b]{30,36;$\pi$}
   \lbl[b]{35,9;$\pi$}
   \lbl[b]{44,36;$v_1$}
   \lbl[b]{70,35;$v_2$}
   \lbl[l]{20,34;$c^*$}
   \lbl[r]{15,40;$n_1$}
   \lbl[r]{10,30;$n_2$}
\end{lpic}
\caption{\label{moje3}}
\end{center}
\end{figure}

We will also need to distinguish between two kinds of surface regions.

\begin{lem}[\cite{Shi98}]
A set $F\backslash\Gamma$ may be \emph{checkerboard coloured}, i.e. we can colour each component with one of two colours $0$ or $1$, such that adjacent (by common boundary edge) regions have different colours.
\end{lem}

We can deduce from this lemma a usefull fact involving our special set of circles.

\begin{cor}
A set $(\bigcup A_f) \cup (F\backslash\Gamma)$ may be checkerboard coloured, because we can colour one by one regions on each component of $F$, treating them as disjont surface knots.
\end{cor}

Let us colour then our set $(\bigcup A_f) \cup (F\backslash\Gamma)$ in (one of two) fixed checkerboard. Finally we will need one more definition to distinguish between two types of elements from the set $A_f$. 

\begin{defn}
Define a set $\emph{X}\subset A_f$ consisting of those circles $d$, that lies on the regions of the colour $1$  and lies higher (with respect to the projection axis) than its mate circle $\left(\pi^{-1}\left(\pi(d)\right)\cap F\right)\backslash d$.
\end{defn}

We are now ready to intruduce our mention invariant.

\begin{thm}\label{tw1}
A class $$\left[\bigcup\left\{c\in A_f:\left(\pi^{-1}\left(\pi(c)\right)\cap F\right)\backslash c\in X\right\}\right]\in \mbox{H}_1\left(F;\mathbb{Z}\right)$$
 is an invariant of $p$-diagrams with respect to the $p$-equivalence of the oriented linked surface $F$.
\end{thm}

\pagebreak

\section{Proofs of theorems}

\begin{proof}[proof of Theorem \ref{tw1}]
We will proceed by checking all four elementary moves that involve $p$-diagrams.

For Roseman moves $I$ and $II$ conclusion of the theorem is satisfied, because both positions and colours of the regions surrounding curves from the set $A_f$ does not change as they projection images do not take part in those moves. In case of the move $III$ we are in the situation of appearing and disappearing of a pair of circles on the double decker set. If those circles do not belong to the set $A_f$ then we give the same argument as in the case of previous moves. But if those circles belong to the set $A_f$ then because their homology class as $1$-chain is zero, our invariant does not change.

More interesting is the last tranformation. In case of the move $IV$, presented in Figure \ref{r43} we can notice that all (involve in this move) curves from double decker set belong either to circles from the set $A_f$ or the are subsets of circles from outside $A_f$ (both before and after triansformation). If it is the latter case, we give argument as in the move $II$ case. Otherwise, we conclude that curves from double decker set (in this local area that we are going to deform) which lies "higher" with respect to the projection direction lies both on the region with the same colour. Moreover this situation does not change after doing this move. Homology class of $1$-chain (with our defined orientation) on which corresponding "lower" curves lies, does not change (because of the group additivity) which ends the proof.
\end{proof}
\begin{proof}[proof of Theorem \ref{t2}]

Let us consider two diagrams $D_F$ and $D_{F'}$ defined by a series of movie stills shown in the Figure \ref{moje2}.

\begin{figure}[ht]
\begin{center}
\begin{lpic}[]{moje2(13cm)}
   \lbl[r]{0,72;$D_{F}$}
   \lbl[r]{0,30;$D_{F'}$}
\end{lpic}
\caption{\label{moje2}}
\end{center}
\end{figure}

They represent $p$-diagrams, because we do not make the Reidemeister third move in any transition between stills. Only such elementary move would produce a triple point in the diagram. Our diagrams give rise to knotted surfaces $F$ and $F'$ that are equivalent. This can be shown for example by applying a combination of another set of moves to the stills. Such \emph{movie moves} are described by Carter and Saito in their book \cite{CarSai98}. However, $p$-diagrams $D_F$ and $D_{F'}$ are not $p$-equivalent. From the Figure \ref{moje1} of a base surface for $F$ with a double decker set drawn, we can conclude that our invariant for the $D_F$ case gives nontrivial class (of one of those circles depicted on the torus) in the oposite to the case of diagram $D_{F'}$ for which the double decker set is empty.
\begin{figure}[ht]
\begin{center}
\begin{lpic}[]{moje1(11cm)}
   \lbl[t]{80,55;$F$}
\end{lpic}
\caption{\label{moje1}}
\end{center}
\end{figure}

\end{proof}

\subsection*{Acknowledgments}

Research of M. Jab\l onowski was partially supported by BW 5107-5-0343-0.

\end{document}